\documentclass[12pt]{amsart}

\usepackage[utf8]{inputenc}
\usepackage[T1]{fontenc}
\usepackage{amsmath,amssymb,amsthm,mathtools}
\usepackage{enumitem}
\usepackage[margin=1.1in]{geometry}
\usepackage[colorlinks=true,linkcolor=blue,citecolor=blue,urlcolor=blue]{hyperref}

\numberwithin{equation}{section}
\allowdisplaybreaks

\newtheorem{theorem}{Theorem}[section]
\newtheorem{proposition}[theorem]{Proposition}
\newtheorem{lemma}[theorem]{Lemma}
\newtheorem{corollary}[theorem]{Corollary}

\theoremstyle{definition}

\theoremstyle{remark}
\newtheorem{remark}[theorem]{Remark}

\newcommand{\im}{\sqrt{-1}}
\newcommand{\C}{\mathbb C}

\newcommand{\PP}{\mathbb P}

\newcommand{\mC}{\mathcal C}
\newcommand{\mO}{\mathcal O}

\newcommand{\Ric}{\operatorname{Ric}}
\newcommand{\HSC}{\operatorname{HSC}}
\newcommand{\rk}{\operatorname{rank}}

\newcommand{\td}{\operatorname{td}}

\newcommand{\tf}{\mathrm{tf}}

\title[HSC and quasi-negative k-Ricci curvature]{Remark on semi-positive holomorphic sectional curvature and quasi-negative $k$-Ricci curvature}

\author{Shiyu Zhang}
\address{School of Mathematical Sciences, University of Science and Technology of China, Hefei, 230026, P.R. China}
\email{shiyu123@mail.ustc.edu.cn}
\date{}

\begin{document}

\begin{abstract}
We record two remarks.  First, for a compact K\"ahler manifold with semi-positive holomorphic sectional curvature, the rational dimension of the MRC fibration is exactly the number of non-truly-flat directions.  Second, for compact K\"ahler manifolds with quasi-negative $k$-Ricci curvature, $1<k<n$, or more generally with quasi-negative mixed curvature $\mC_{a,b}$ for $a,b>0$, the canonical bundle is ample. 
\end{abstract}
\maketitle

\tableofcontents

\section{Introduction}
Throughout this note, $(X,g)$ denotes a compact K\"ahler manifold of complex dimension $n$, and $R$ denotes the Chern curvature tensor of $g$.

\subsection{Semi-positive holomorphic sectional curvature}
The problem for rational connectedness of compact K\"ahler manifolds with positive holomorphic sectional curvature ($\HSC>0$) goes back to Yau's problem \cite[Problem 47]{Yau82}. It was completely resolved in \cite{yang2018}. Yau's conjecture also holds in the quasi-positive case \cite{heier2020,matsumura2022} for projective manfiolds, which was extended to the K\"ahler setting in \cite{ZZ25}. Recently, the structure theorems of compact K\"ahler manifolds  with $\HSC\geq 0$ in terms of MRC fibrations has been extablished in \cite{Matsumura2025}, we refer to \cite{MatsumuraSurvey} for a detailed introduction on this topic.

For $x\in X$, a tangent vector $v\in T^{1,0}_xX$ is called \emph{truly flat} if
\[
        R(v,\overline u,y,\overline z)=0,
        \qquad \forall u,y,z\in T^{1,0}_xX.
\]
Set
\[
        V_{\tf,x}:=\{v\in T^{1,0}_xX\mid v\text{ is truly flat}\},
        \qquad
        n_{\tf}(X,g):=n-\inf_{x\in X}\dim_{\C}V_{\tf,x}.
\]
The invariant $n_{\tf}(X,g)$ measures the number of non-flat directions of the curvature tensor. Matsumura (\cite[Theorem 1.1]{Matsumura2025}) obtained the structure theorems for a compact K\"ahler manifold with $\HSC\geq 0$ and derived the following inequality
\[
        \dim X-\dim Y\ge n_{\tf}(X,g)
\]
for an MRC fibration $f:X\dashrightarrow Y$; see \cite{Matsumura2025}.  The remaining equality is the missing part of the structure problem for compact K\"ahler manifolds with non-negative holomorphic sectional curvature; see \cite[Remark~3.7]{MatsumuraSurvey}.

The first observation of this note is the reverse inequality (Proposition \ref{prop:rational-dim-upper}). Then we obtain

\begin{theorem}\label{thm:rational-dim-equality}
Let $(X,g)$ be a compact K\"ahler manifold with $\HSC\geq0$. Then for any MRC fibration $f:X\dashrightarrow Y$,
$$\dim X-\dim Y=n_{\tf}(X,g).$$
\end{theorem}

This is the optimal version of Yau's conjecture on $\HSC\geq 0$ and can be viewed as the counterpart of \cite[Conjecture 1.11]{heier2018}, which predicts an optimal version of Wu-Yau's theorem on $\HSC\leq 0$, i.e., 
$$\mathrm{Kod}(X)=n_{\tf}(X,g).$$

Combining with \cite[Theorem 1.1]{Matsumura2025}, we obtain
\begin{corollary}\label{coro:splitting}
       Let $(X,g)$ be a compact K\"ahler manifold with $\HSC\geq0$. Then the universal cover $(\widehat{X},\widehat{g})$ of $(X,g)$ splits holomorphically and isometrically as
$$(\widehat{X},\widehat{g})\cong (\C^m,g_{\C^m})\times (F,g_F)$$
for some rationally connected manifold $(F,g_F)$ with $\HSC_{g_F}\geq0$ and $\dim F=n_{\tf}(X,g)$;
\end{corollary}

$\dim F=n_{\tf}(X,g)$ is a quasi-positive curvature condition, when $\HSC\geq 0$, it is equivalent to that $(T_F,g_F)$ is RC-quasi-positive (cf. Remark \ref{rem:quasi-positive}).

\subsection{Quasi-negative \texorpdfstring{$k$}{k}-Ricci curvature and ampleness}
Let $\Sigma\subset T^{1,0}_{X,x}$ be a $k$-dimensional holomorphic subspace and let $\{e_i\}_{i=1}^k$ be a unitary frame of $\Sigma$.  The $k$-Ricci curvature along $\Sigma$ is
\[
        \Ric_k(x,\Sigma)(v,\overline v)
        :=\sum_{i=1}^k R(e_i,\overline{e_i},v,\overline v),
        \qquad v\in \Sigma.
\]
The $k$-Ricci curvature was introduced by Ni in the study of $k$-hyperbolicity of compact K\"ahler manifolds \cite{nicpam}.  It coincides with holomorphic sectional curvature when $k=1$ and with Ricci curvature when $k=n$.  We say that $\Ric_k$ is non-positive if this Hermitian form is non-positive for every $x$ and every $\Sigma$, and quasi-negative if it is non-positive everywhere and negative at least at one point in the usual sense.

For $a,b>0$ and $\phi\in C^{\infty}(X)$, define the mixed curvature
\[
        \mC_{a,b,\phi}(v)
        :=a\Ric_g(v,\overline v)
        +b\frac{R(v,\overline v,v,\overline v)}{|v|_g^2}
        +\im\partial\overline\partial\phi(v,\overline v),
        \qquad v\ne0.
\]
We write $\mC_{a,b}\le0$ if such a $\phi$ exists and $\mC_{a,b,\phi}(v)\le0$ for every $v\ne0$.  The quasi-negative condition is defined similarly.

For $k=1$, Wu--Yau \cite{WY1} and Tosatti--Yang \cite{TY17}, proved the ampleness of the canonical bundle under negative holomorphic sectional curvature.  More generally, Diverio--Trapani \cite{DT19} for $k=1$ and Chu--Lee--Tam \cite{CLT22} for $1<k<n$ proved that quasi-negative $k$-Ricci curvature implies that $K_X$ is nef and big; compare also \cite{WY2,LNZ21}.  Royden's observation \cite{Royden80} shows that non-positive holomorphic sectional curvature rules out rational curves, so the ampleness in the case $k=1$ follows from Mori's cone theorem.  For $1<k<n$, the absence of rational curves is not known a priori, and Chu--Lee--Tam imposed it as an additional hypothesis to conclude the ampleness \cite[Theorem 1.1]{CLT22}.  The second observation removes this hypothesis.

\begin{theorem}\label{thm:main-kricci-intro}
Let $(X,g)$ be a compact K\"ahler manifold. If $\mC_{a,b}$ is quasi-negative for some $a,b>0$ (which is satisfied when $\Ric_k$ is quasi-negative for some $1<k<n$), then $K_X$ is ample.
\end{theorem}

The key point is that non-positive mixed curvature forces a strictly positive lower bound for $K_X\cdot C$ on every rational curve $C\subset X$.  Combining this with the bigness theorem of Chu--Lee--Tam and the cone theorem gives ampleness.

\vspace{0.6cm}

\noindent\textbf{AI disclosure.}
   The proof of Proposition \ref{prop:rational-dim-upper} is motivated by GPT-5.5 Pro.

\section{Rational dimension and truly flat tangent vectors}

We prove the reverse inequality needed for Theorem~\ref{thm:rational-dim-equality}.  The argument is independent of the sign of curvature.

\begin{proposition}\label{prop:rational-dim-upper}
Let $X$ be a compact K\"ahler manifold and let $h$ be any Hermitian metric on $T_X$.  For any MRC fibration $f:X\dashrightarrow Y$,
\[
        \dim X-\dim Y\le n_{\tf}(X,h).
\]
\end{proposition}

We first record an elementary observation about Todd forms.

\begin{lemma}\label{lem:todd-vanishing}
Let $(E,h)$ be a Hermitian vector bundle over a complex manifold $X$ of dimension $n$.  For $x\in X$, set
\[
        n_{\tf}(E,h,x)
        := n-\rk\{v\in T^{1,0}_{X,x}\mid \Theta^h(v,\overline u)=0,
        \ \forall u\in T^{1,0}_{X,x}\},
\]
where $\Theta^h$ is the Chern curvature tensor of $h$.  If $s>n_{\tf}(E,h,x)$, then
\[
        \td_s(E,h)(x)=0.
\]
\end{lemma}

\begin{proof}
The Todd form is
\[
        \td(E,h)=\det\left(\frac{\im\Theta^h}{1-e^{-\im\Theta^h}}\right)
        =1+\sum_{s=1}^{n}\td_s(E,h).
\]
Only the directions on which $\Theta^h$ is nonzero can contribute to positive-degree components.  Hence the component of degree $s$ vanishes at $x$ whenever $s$ is larger than the rank of the non-flat part.
\end{proof}

\begin{proof}[Proof of Proposition~\ref{prop:rational-dim-upper}]
Set $d:=n_{\tf}(X,h)$ and argue by contradiction.  Suppose
\[
        k:=\dim X-\dim Y>d.
\]
By Lemma~\ref{lem:todd-vanishing},
\[
        \td_k(T_X,h)=0
\]
as a differential form, and hence
\[
        \td_k(T_X)=0\qquad\text{in }H^{2k}(X,\C).
\]
Take a smooth model of the MRC fibration and a general smooth fiber $F$ of dimension $k$.  Over a nonempty open set of the base, the fibration is a holomorphic submersion, and restricting
\[
        0\longrightarrow T_{X/Y}\longrightarrow T_X\longrightarrow f^*T_Y\longrightarrow0
\]
to $F$ gives
\[
        0\longrightarrow T_F\longrightarrow T_X|_F\longrightarrow \mO_F^{\oplus \dim Y}\longrightarrow0.
\]
Thus
\[
        i^*\td(T_X)=\td(T_X|_F)=\td(T_F),
\]
where $i:F\hookrightarrow X$ is the inclusion.  Hirzebruch--Riemann--Roch formula gives
\[
        \chi(F,\mO_F)
        =\int_F\td(T_F)
        =\int_F\td_k(T_F)
        =\int_F i^*\td_k(T_X)=0.
\]
However, a general fiber of an MRC fibration is rationally connected, so
\[
        \chi(F,\mO_F)=1,
\]
a contradiction.
\end{proof}

\begin{proof}[Proof of Theorem~\ref{thm:rational-dim-equality} and Corollary \ref{coro:splitting}]
Theorem~\ref{thm:rational-dim-equality} can be deduced by Proposition \ref{prop:rational-dim-upper} and \cite[Corollary 1.3]{Matsumura2025}. Corollary \ref{coro:splitting} can be implied by Theorem \ref{thm:rational-dim-equality} and \cite[Theorem 1.1]{Matsumura2025}, since $\phi$ in the statement of \cite[Theorem 1.1]{Matsumura2025} is an MRC fibration.
\end{proof}

\begin{remark}[RC-quasi-positivity]\label{rem:quasi-positive}
    Let $(T_F,g_F)$ be a compact K\"ahler manifold with $\HSC_{g_F}\geq 0$ and $\dim F=n_{\tf}(F,g_F)$. Suppose that $(T_F,g_F)$ is not RC-quasi-positive, then for any $x\in F$, there exists some $v\in T_{F,x}\setminus {0}$ such that
    $$R(v,\bar{v},u,\bar{u})\leq 0,\ \forall u\in T_{F,x}.$$
    In particular, $\HSC(v)=0$ and thus \cite[Lemma 6.1]{yang2018} implies
    $$R(v,\bar{v},u,\bar{u})\geq 0,\ \forall u\in T_{F,x}.$$
    Therefore,
    $$R(v,\bar{v},u,\bar{u})=0,\ \forall u\in T_{F,x},$$
    which implies $v\in V_{\tf,x}$ and concludes a contradiction.
\end{remark}

\section{Quasi-negative \texorpdfstring{$k$}{k}-Ricci curvature}

We now prove Theorem \ref{thm:main-kricci-intro}.

\subsection{Mixed curvature and rational curves}
Recall the following comparison from Chu--Lee--Tam.

\begin{lemma}\label{lem:kricci-to-mixed}
If $\Ric_k$ is non-positive, respectively quasi-negative, then $\mC_{k-1,n-k}$ is non-positive, respectively quasi-negative.
\end{lemma}

\begin{proof}
This is \cite[Lemma~2.2]{CLT22}.
\end{proof}

The key point is the following standard estimate.

\begin{proposition}\label{prop:rational-curve-bound}
Let $(X,g)$ be a compact K\"ahler manifold with non-positive mixed curvature $\mC_{a,b}\le0$ for some $a,b>0$.  Then every rational curve $C\subset X$ satisfies
\[
        K_X\cdot C\ge \frac{2b}{a}.
\]
In particular, $K_X\cdot C>0$ for every rational curve $C$.
\end{proposition}

The proof built on the Gauss-Codazzi equation.

\begin{proof}
Let $f:\PP^1\to X$ be a nonconstant holomorphic map whose image is $C$.  The differential gives a morphism
\[
        df:T_{\PP^1}\longrightarrow f^{-1}T_X.
\]
Let $L$ be the saturation of the image of $df$.  Then $L$ is a line subbundle of $f^{-1}T_X$ equipped with the induced metric $h_L$.  On the dense open set $W\subset\PP^1$ where $df\ne0$, the bundle $L$ is identified with $T_{\PP^1}$ and $h_L=f^*g$.  Therefore
\[
        2=\deg T_{\PP^1}\le \deg L
        =\int_{\PP^1}c_1(L,h_L)
        =\int_W S(f^*g)\,f^*\omega_g,
\]
where $S(f^*g)$ denotes the scalar curvature of the induced metric on the curve.

Fix a point of $f(W)$ and a nonzero tangent vector $v$ tangent to the image curve.  The curvature assumption gives
\[
        b\frac{R(v,\overline v,v,\overline v)}{|v|_g^2}
        +a\Ric_g(v,\overline v)
        +\im\partial\overline\partial\phi(v,\overline v)
        \le0.
\]
By the Gauss--Codazzi equation, the intrinsic curvature of the curve is bounded above by the ambient holomorphic sectional curvature in the tangent direction.  Hence
\begin{align*}
        \int_W S(f^*g)\,f^*\omega_g
        &\le \int_{f(W)}\frac{R(v,\overline v,v,\overline v)}{|v|_g^4}\,\omega_g|_{f(W)}\\
        &\le \frac{a}{b}\int_{f(W)}(-\Ric_g)|_{f(W)}
        -\frac1b\int_{f(W)}\im\partial\overline\partial\phi|_{f(W)}\\
        &=\frac{a}{b}\,K_X\cdot C.
\end{align*}
The last term involving $\partial\overline\partial\phi$ integrates to zero on $\PP^1$.  Combining with the previous inequality gives
\[
        2\le \frac{a}{b}K_X\cdot C,
\]
so $K_X\cdot C\ge \frac{2b}{a}$.
\end{proof}

\begin{remark}
Chu--Lee--Tam proved nefness of $K_X$ for $\mC_{a,b}\le0$ using the twisted K\"ahler--Ricci flow \cite{CLT22}.  Proposition~\ref{prop:rational-curve-bound} gives a direct rational-curve obstruction to $K_X$ being negative on an extremal rational curve.
\end{remark}

\subsection{Ampleness}
We use the following standard consequence of the cone theorem and the base point free theorem.

\begin{lemma}\label{lem:big-positive-rational-curves-ample}
Let $X$ be a projective manifold of general type.  If
\[
        K_X\cdot C>0
\]
for every rational curve $C\subset X$, then $K_X$ is ample.
\end{lemma}

\begin{proof}
If $K_X$ is not nef, Mori's cone theorem gives a rational curve $C$ with $K_X\cdot C<0$ \cite[Theorem~7.49]{De2001}, a contradiction.  Hence $K_X$ is nef.  Since $X$ is of general type, $K_X$ is nef and big, so it is semiample by the base point free theorem \cite[Theorem~7.32]{De2001}.

If $K_X$ is not ample, then there is a curve $C'$ with $K_X\cdot C'=0$.  Equivalently, for a suitable effective divisor $D$ and every sufficiently small rational $\epsilon>0$,
\[
        (K_X+\epsilon D)\cdot C'<0,
\]
with $(X,\epsilon D)$ klt; see, for example, \cite[Lemma~2.1]{DT19}.  The logarithmic cone theorem \cite[Theorem~7.49]{De2001} gives a decomposition of the class of $C'$ modulo the nonnegative part of the cone involving finitely many rational curves $C_i$ with positive coefficients.  Since $K_X$ is nef and $K_X\cdot C_i>0$ for all rational curves $C_i$, this forces
\[
        K_X\cdot C'>0,
\]
contradicting $K_X\cdot C'=0$.
\end{proof}

\begin{theorem}\label{thm:mixed-ample}
Let $(X,g)$ be a compact K\"ahler manifold with quasi-negative $\mC_{a,b}$ for some $a,b>0$.  Then $K_X$ is ample.
\end{theorem}

\begin{proof}
Chu--Lee--Tam proved that quasi-negative $\mC_{a,b}$ implies that $K_X$ is big \cite[Theorem~1.1(3)]{CLT22}.  Since $X$ is compact K\"ahler and carries a big line bundle, $X$ is Moishezon; hence $X$ is projective.  Proposition~\ref{prop:rational-curve-bound} shows that $K_X\cdot C>0$ for every rational curve $C\subset X$.  Lemma~\ref{lem:big-positive-rational-curves-ample} then gives that $K_X$ is ample.
\end{proof}

\begin{proof}[Proof of Theorem~\ref{thm:main-kricci-intro}]
The mixed-curvature statement is Theorem~\ref{thm:mixed-ample}.  If $\Ric_k$ is quasi-negative for $1<k<n$, then Lemma~\ref{lem:kricci-to-mixed} gives quasi-negative $\mC_{k-1,n-k}$, where both coefficients are positive.  Applying Theorem~\ref{thm:mixed-ample} proves that $K_X$ is ample.
\end{proof}

\end{document}